\pdfoutput=1
\RequirePackage{ifpdf}
\ifpdf 
\documentclass[pdftex]{sigma}
\else
\documentclass{sigma}
\fi

\numberwithin{equation}{section}

\newtheorem{Theorem}{Theorem}[section]
\newtheorem{Corollary}[Theorem]{Corollary}
\newtheorem{Lemma}[Theorem]{Lemma}
\newtheorem{Proposition}[Theorem]{Proposition}
 { \theoremstyle{definition}
\newtheorem{Definition}[Theorem]{Definition}
\newtheorem{Construction}[Theorem]{Construction}
\newtheorem{Example}[Theorem]{Example}
\newtheorem{Remark}[Theorem]{Remark} }

\usepackage{tikz}
\usetikzlibrary{arrows,chains,matrix,positioning,scopes, decorations.markings}

\newcommand {\TT} {{\T\in\mathcal{T}}}

\newcommand {\on}[1] {\operatorname{#1}}
\newcommand {\pt} {\on{pt}}
\newcommand {\mmod} {/\!\!/}
\newcommand {\ZZ} {\mathbb Z}
\renewcommand {\(} {\left(}
\renewcommand {\)} {\right)}
\newcommand {\tensor}{\otimes}
\newcommand {\CC} {\mathbb C}
\newcommand {\eps} {\varepsilon}
\newcommand {\id} {\operatorname{id}}
\newcommand {\Hom} {\operatorname {Hom}}
\newcommand {\isomap} [3] {\specialmap {#1} {#2}{\overset {\cong{\phantom{.}}}
 {\longrightarrow}} {#3}}
\newcommand {\specialmap} [4] {\text {$ #1\colon #2 #3 #4 $}}
\newcommand {\map} [3] {\specialmap {#1} {#2}{\to} {#3}}
\newcommand {\sub} {\subseteq}
\newcommand {\longmap} [3] {\specialmap {#1} {#2}{\longrightarrow} {#3}}

\newcommand {\inv} {^{-1}}

\newcommand {\LE} {\mathbf L}

\newcommand {\mC} {\mathcal C}

\newcommand {\mG} {\mathcal G}

\newcommand {\mI} {\mathcal I}
\newcommand {\mL} {\mathcal L}

\newcommand {\mP} {\mathcal P}

\newcommand {\mT} {\mathcal T}
\newcommand {\mU} {\mathcal U}

\newcommand {\ul}[1]{\underline{#1}}

\newcommand {\bbS} {{\mathbb S}}
\newcommand {\FF} {{\mathbb F}}

\newcommand {\RR} {{\mathbb R}}

\renewcommand {\TT} {\mathbb T} 

\newcommand {\vL} {{\cech{\Lambda\!}}}
\newcommand {\Lat} {{\vL}}

\usepackage [all]{xy}
\SelectTips{cm}{12}
\usepackage[T1]{fontenc}

\def\cech#1{#1^{\smash{\scalebox{.7}[1.4]{\rotatebox{90}{\guilsinglleft}}}}}

\newcommand{\vC}
 {{C\hspace{-1.205ex}^{\raisebox{2.9pt}{\smash{\scalebox{.9}[.9]{\rotatebox{90}
 {\guilsinglleft}}}}}\hspace{.4pt}}}
\newcommand{\vH}
 {{H\hspace{-1.32ex}^{\raisebox{2.55pt}{\smash{\scalebox{.9}[.9]{\rotatebox{90}
 {\guilsinglleft}}}}}} \hspace{1pt}}

\newcommand{\vT}
 {{T\hspace{-1.2415ex}^{\raisebox{2.55pt}{\smash{\scalebox{.9}[.9]{\rotatebox{90}{\guilsinglleft}}}}}}}

\pgfdeclareshape{hackgenuspic}{
 \anchor{center}{\pgfpointorigin}
\backgroundpath{
 \pgfpathmoveto{\pgfqpoint{-0.78cm}{-.17cm}}
 \pgfpathcurveto %
 {\pgfpoint{-0.35cm}{-.44cm}}
 {\pgfpoint{0.35cm}{-.44cm}}
 {\pgfpoint{.78cm}{-0.17cm}}
 \pgfpathmoveto{\pgfqpoint{-0.78cm}{-0.17cm}}
 \pgfpathcurveto %
 {\pgfpoint{-0.25cm}{.25cm}}
 {\pgfpoint{.25cm}{.25cm}}
 {\pgfpoint{0.78cm}{-0.17cm}}
 \pgfsetfillcolor{white}
 \pgfusepath{fill}

 \pgfpathmoveto{\pgfqpoint{-1cm}{0cm}}
 \pgfpathcurveto %
 {\pgfpoint{-0.5cm}{-.5cm}}
 {\pgfpoint{0.5cm}{-.5cm}}
 {\pgfpoint{1cm}{0cm}}

 \pgfpathmoveto{\pgfqpoint{-0.75cm}{-0.15cm}}
 \pgfpathcurveto %
 {\pgfpoint{-0.25cm}{.25cm}}
 {\pgfpoint{.25cm}{.25cm}}
 {\pgfpoint{0.75cm}{-0.15cm}}
 \pgfusepath{stroke}
 }
 }

\tikzstyle directed=[postaction={decorate,decoration={markings,
 mark=at position .65 with {\arrow{stealth}}}}]
\tikzstyle reverse directed=[postaction={decorate,decoration={markings,
 mark=at position .65 with {\arrowreversed{stealth};}}}]

\makeatletter
\tikzset{join/.code=\tikzset{after node path={%
\ifx\tikzchainprevious\pgfutil@empty\else(\tikzchainprevious)%
edge[every join]#1(\tikzchaincurrent)\fi}}}
\makeatother
\tikzset{>=stealth',every on chain/.append style={join},
 every join/.style={->}}
\tikzstyle{labeled}=[execute at begin node=$\scriptstyle,
 execute at end node=$]

\begin{document}
\allowdisplaybreaks

\newcommand{\arXivNumber}{1406.7046}

\renewcommand{\thefootnote}{}

\renewcommand{\PaperNumber}{014}

\FirstPageHeading

\ShortArticleName{Categorical Tori}

\ArticleName{Categorical Tori\footnote{This paper is a~contribution to the Special Issue on Moonshine and String Theory. The full collection is available at \href{https://www.emis.de/journals/SIGMA/moonshine.html}{https://www.emis.de/journals/SIGMA/moonshine.html}}}

\Author{Nora GANTER}

\AuthorNameForHeading{N.~Ganter}

\Address{School of Mathematics and Statistics, The University of Melbourne,\\ Parkville, Victoria 3010, Australia}
\Email{\href{mailto:noraganter@gmail.com}{noraganter@gmail.com}}
\URLaddress{\url{http://researchers.ms.unimelb.edu.au/~nganter@unimelb/}}

\ArticleDates{Received September 23, 2017, in f\/inal form January 31, 2018; Published online February 17, 2018}

\Abstract{We give explicit and elementary constructions of the categorical extensions of a~torus by the circle and discuss an application to loop
group extensions. Examples include maximal tori of simple and simply connected compact Lie groups and the tori associated to the Leech and Niemeyer lattices. We obtain the extraspecial 2-groups as the isomorphism classes of categorical f\/ixed points under an involution action.}


\Classification{22E99; 18D99}

\renewcommand{\thefootnote}{\arabic{footnote}}
\setcounter{footnote}{0}

\section{Introduction}
By a {\em categorical group}, or a {\em $2$-group}, we mean a monoidal groupoid $(\mG,\bullet,1)$ with weakly invertible objects. If $\mG$ is a Lie groupoid, the monoidal structure is required to be locally smooth in an appropriate sense \cite{Schommer-Pries:Central_extensions}. In this situation, one speaks of a {\em Lie $2$-group}. Categorical groups play as important a~role in string theory as groups do in particle physics \cite{Baez:Huerta:An_invitation}, and
a number of prominent groups seem to be most naturally understood via their categorical ref\/inements. Most famously, the inf\/inite-dimensional groups $\on{String(n)}$ come from f\/inite-dimensional Lie 2-groups \cite{Schommer-Pries:Central_extensions}. Also, Weyl groups and some of the sporadic
groups, including the Monster, are known or conjectured to be the isomorphism classes of categorical groups (see Section~\ref{sec:examples}). In recent years, there is a~growing understanding that the categorical nature of these groups is worth exploring, and we will see that this perspective sheds new light on some old and important mathematics. The purpose of this paper is to give a simple and hands-on description of a basic class of examples, namely all central extensions of Lie 2-groups of the form
 \begin{equation*}
 \begin{tikzpicture}
 \node at (0,0) [name=a] {$1$};
 \node at (2,0) [name=b] {$\pt\mmod U(1)$};
 \node at (4.4,0) [name=c] {$\mT$};
 \node at (6.3,0) [name=d] {$T$};
 \node at (7.75,0) [name=e] {$1$,};
 \draw[->] (a) -- (b);
 \draw[->] (b) -- (c);
 \draw[->] (c) -- (d);
 \draw[->] (d) -- node [midway, above] {} (e);
\end{tikzpicture}
\end{equation*}
where $T$ is a compact torus, and $\pt\mmod U(1)$ is the one-object groupoid with $\on{Aut}(\pt)=U(1)$. We will refer to such a $\mathcal T$ as a {\em categorical torus}. Categorical tori are important for a number of reasons: f\/irst, they turn up as maximal tori of Lie 2-groups. So, any character theory of Lie 2-groups should sensibly start here. Second, categorical tori incorporate many aspects relevant to the construction of
sporadic groups (see Sections~\ref{sec:examples} and~\ref{sec:Extraspecial}). Finally, the theory of Lie 2-groups is closely related to that
of loop groups, a point we will explore in Section~\ref{sec:loops}. Our approach is to work from 2-groups to loop groups. In particular, our f\/irst construction of $\mT$ does not involve loops or any other inf\/inite-dimensional objects. The {transgression} machine of~\cite{Waldorf:Transgression_to_loop_spaces_and_its_inverse_II} then yields a~simple description of central extensions of the loop group~$\mL T$. What we hope to get across is that categorical groups can be fairly simple to construct and easy to work with and that some of the complicated features of loop groups are merely the shadow of rather obvious phenomena on the categorical side. This is the f\/irst in a series of papers describing the representation theory of categorical tori.

\section{Constructions of categorical tori}\label{sec:cons}
Let $\Lat$ be a free $\ZZ$-module of f\/inite rank, and let $J$ be an integer-valued bilinear form on $\Lat$. Tensoring $\Lat$ over $\ZZ$ with the short exact sequence
\begin{equation*}
\begin{tikzpicture}
 \node at (0.2,0) [name=a] {$0$};
 \node at (2,0) [name=b] {$\ZZ$};
 \node at (4,0) [name=c] {$\RR$};
 \node at (6.2,0) [name=d] {$U(1)$};
 \node at (8.2,0) [name=e] {$0$,};
 \node at (4,-1) [name=f] {$r$};
 \node at (6.2,-1) [name=g] {$e^{2\pi ir}$,};
 \draw[->] (a) -- (b);
 \draw[->] (b) -- (c);
 \draw[->] (c) -- node [midway, above] {$\exp$} (d);
 \draw[->] (d) -- (e);
 \draw[|->] (f) -- (g);
\end{tikzpicture}
\end{equation*}
we obtain the short exact sequence
\begin{equation*}
\begin{tikzpicture}
 \node at (0.2,0) [name=a] {$0$};
 \node at (2,0) [name=b] {$\Lat$};
 \node at (4,0) [name=c] {$\mathfrak t$};
 \node at (6.2,0) [name=d] {$T$};
 \node at (8.2,0) [name=e] {$0$.};
 \draw[->] (a) -- (b);
 \draw[->] (b) -- (c);
 \draw[->] (c) -- node [midway, above] {$\exp$} (d);
 \draw[->] (d) -- (e);
\end{tikzpicture}
\end{equation*}
Here $T$ is the compact torus with coweight lattice $\Lat$, and $\mathfrak t$ is its Lie-algebra. We will give three equivalent constructions of the categorical torus associated to $(\Lat,J)$.
\begin{Construction}\label{cons:strict}
 Let $\mathfrak t$ act on $\Lat\times U(1)$ by
 \begin{equation*}
 \begin{tikzpicture}
 \node at (0,0) [name=A, anchor=east] {$\mathfrak t\times\Lat\times U(1)$};
 \node at (1,0) [name=B, anchor=west] {$\Lat\times U(1)$,};
 \node at (0,-0.8) [name=C, anchor=east] {$(x,m,z)$};
 \node at (1,-0.8) [name=D, anchor=west] {$(m,z\cdot\exp(J(m,x)))$.};
 \draw[->] (A) -- (B);
 \draw[|->] (C) -- (D);
 \end{tikzpicture}
 \end{equation*}
 Then the strict categorical group $(\mT,\bullet, 0)$ associated
 to the crossed module
 \begin{equation*}
 \begin{tikzpicture}
 \node at (0,0) [name=A, anchor=east] {\ensuremath{\Psi\colon \ \Lat\times U(1)}};
 \node at (1,0) [name=B, anchor=west] {\ensuremath{\mathfrak t,}};
 \node at (0,-0.8) [name=C, anchor=east] {\ensuremath{(m,z)}};
 \node at (1,-0.8) [name=D, anchor=west] {\ensuremath{m}};
 \draw[->] (A) -- (B);
 \draw[|->] (C) -- (D);
 \end{tikzpicture}
 \end{equation*}
 is a categorical torus.
\end{Construction}
Explicitly, the Lie groupoid $\mT$ has objects $\mathfrak t$ and
arrows $\mathfrak t\ltimes(\Lat\times U(1))$, which we write as
\begin{equation*}
 \begin{tikzpicture}
 \node at (0,0) [name=a] {$x$};
 \node at (1.8,0) [name=b] {$x+m$};
 \draw[->] (a) -- node [midway, above] {$z$} (b);
 \end{tikzpicture}
\end{equation*}
with $x\in \mathfrak t$, $m\in\Lat$ and $z\in U(1)$, to indicate
source and target. Composition in $\mT$ is
\begin{equation*}
 \begin{tikzpicture}
 \node at (0,0) [name=a] {$x$};
 \node at (1.8,0) [name=b] {$x+m$};
 \node at (4.3,0) [name=c] {$x+m+n$,};
 \draw[->] (a) -- node [midway, above] {$z$} (b);
 \draw[->] (b) -- node [midway, above] {$w$} (c);
 \draw[->] (a) .. controls (1.5,-1) and (2.8,-1) .. node [midway, below] {$zw$} (c);
 \end{tikzpicture}
\end{equation*}
and tensor multiplication is
\begin{gather*}
 x\bullet y = x+y, \\
 \left(x\xrightarrow{\,\,\,\,z\,\,\,\,}x+m\right)\bullet
 \left(y\xrightarrow{\,\,\,\,w\,\,\,\,}y+n\right) =
 \left(x+y\xrightarrow{\,\,\,\,zw\exp(J(m,y))\,\,\,\,}x+y+m+n\right).
\end{gather*}
The unit object is $0$, and the associativity and unit isomorphisms are identities. This gives $\mT$ the structure of a group object in the
category of Lie groupoids. When we want to emphasize the dependence on $J$, we will use the notation $\bullet_J$ for $\bullet$. In the following, we will use the notation~$X\mmod G$ for the translation groupoid associated to the action of a group $G$ on a space $X$. With this notation, the underlying groupoid of Construction \ref{cons:strict} is
\begin{gather*}
 \mT \cong (\mathfrak t\mmod\Lat )\times (\pt\mmod U(1) ),
\end{gather*}
where $\Lat$ acts on $\mathfrak t$ by translations.
\begin{Construction}\label{cons:compact}
 Consider the torus $T$ as a Lie groupoid with only identity arrows, and let
 \begin{gather*}
 p\colon \ \mathfrak t\mmod\Lat \xrightarrow{\,\,\sim\,\,} T
 \end{gather*}
 be the equivalence of Lie-groupoids sending the object $x$ to $\exp(x)$. Note that $p$ does not possess a continuous inverse. The language developed in \cite{Schommer-Pries:Central_extensions} interprets \begin{gather*}
 \mathfrak T = T\times \pt\mmod U(1),
 \end{gather*}
 together with the data inherited from
 Construction \ref{cons:strict}, as a 2-group object with multiplication
 \begin{equation} \label{eq:zig-zag}
 \mathfrak T\times\mathfrak T \xleftarrow{\,\,\,\,\sim\,\,\,\,}
 \mT\times\mT \xrightarrow{\,\,\,\,\bullet\,\,\,\,} \mT
 \xrightarrow{\,\,\,\,\sim\,\,\,\,} \mathfrak T
 \end{equation}
 in a suitable
 localization of the bicategory of Lie-groupoids. The equivalences in \eqref{eq:zig-zag} are those induced by $p$. Dif\/ferent communities have dif\/ferent language for the 1-morphisms in this localized bicategory $\mathbf{Bibun}$. Depending on taste, the reader may wish to think of this multiplication on~$\mathfrak T$ as a {\em zig-zag}, a {\em span}, an {\em orbifold map}, an {\em anafunctor}, or a {\em bibundle}.
\end{Construction}
Our third construction is as a multiplicative bundle gerbe in the sense of \cite{Brylinski:Loop_spaces_characteristic_classes} and~\cite{Carey:Johnson:Murray:Stevenson:Wang}.
\begin{Construction}\label{cons:mbg}
 Let $\LE$ be the line bundle over $T\times T$ with multipliers
 \begin{align*}
 f_{(m,n)}\colon \ \mathfrak t\times\mathfrak t & \longrightarrow U(1),\\
 (x,y) & \longmapsto \exp\({J(m,y)}\),
 \end{align*}
 $(m,n)\in\Lat\times\Lat$.
 We claim that over $T\times T\times T$ we have a canonical isomorphism
 \begin{gather*}
 \alpha\colon \ \on{m}_{12}^*\LE\tensor \on{pr}_{12}^*\LE \cong
 \on{m}_{23}^*\LE\tensor \on{pr}_{23}^*\LE,
 \end{gather*}
 where $\on{pr}_{ij}$ is the projection onto factors $i$ and $j$ and $\on{m}_{ij}$ is given by multiplication of these two factors. Indeed, source and target of $\alpha$ have identical multipliers
 \begin{align*}
 f_{(k,m,n)}\colon \ \mathfrak t\times\mathfrak t\times\mathfrak t & \longrightarrow U(1),\\
 (x,y,z) & \longmapsto \exp(J(k,y) + J(k,z) + J(m,z)),
 \end{align*}
 $(k,m,n)\in\Lat\times\Lat\times\Lat$. The pair $(\LE,\alpha)$ equips the trivial bundle gerbe over $T$ with a multiplicative structure.
 It is well known (e.g., \cite[Theorem~3.2.5]{Waldorf:A_construction_of_string_2-group_models}) that the data of a multiplicative bundle gerbe over $T$ are equivalent to those of a Lie 2-group extension of $T$ by $\pt\mmod U(1)$.
\end{Construction}
Explicitly, the line bundle $\LE$ is constructed as
 \begin{gather*}
 \LE = \CC\times \mathfrak t\times\mathfrak t /{\sim}
 \end{gather*}
 with
 \begin{gather*}
 \(z,x,y\) \sim (z\cdot \exp(J(m,y)),x+m,y+n)
 \end{gather*}
 for $(m,n)\in \Lat\times\Lat$. The trivial bundle gerbe over $T$ corresponds to the groupoid $\mI$ whose objects are pairs $(t,L)$ with $t\in T$
and $L$ a hermitian line, and whose arrows \begin{gather*}(t,L)\longrightarrow (t,L')\end{gather*}
are the unitary isomorphisms from $L$ to $L'$. The line bundle $\LE$ equips $\mI$ with the monoidal structure
\begin{gather*}
 (s,L_1)\circledast(t,L_2) = (s\cdot t,\mathbf L_{s,t}\tensor L_1\tensor L_2),
\end{gather*}
whose associativity isomorphisms
\begin{gather*}
 \alpha_{r,s,t}\colon \ \mathbf L_{r\cdot s,t}\tensor \mathbf L_{r,s} \cong \mathbf
 L_{r,s\cdot t}\tensor \mathbf L_{s,t}
\end{gather*}
are encoded in $\alpha$.
\begin{Proposition}
 The three constructions yield equivalent extensions of $T$ by $\pt\mmod U(1)$.
\end{Proposition}
\begin{proof}
It is clear that Construction \ref{cons:strict} and Construction~\ref{cons:compact} are equivalent. To see their equivalence with Construction~\ref{cons:mbg}, let $F$ be the functor
 \begin{align*}
 F\colon \ \mT & \longrightarrow \mI,\\
 x & \longmapsto (\exp(x),\CC) \qquad \text{on objects},\\
 (x,m,z) & \longmapsto z \hspace{24.5mm} \text{on arrows}.
 \end{align*}
 We want to make $F$ into a monoidal equivalence. For this, we need to specify an isomorphism
 \begin{gather*}
 \eps\colon \ (1,\CC) \cong F(0)
 \end{gather*}
 and a natural isomorphism
 \begin{gather*}
 \phi\colon \ F(-)\circledast F(-) \Longrightarrow F(-\bullet-),
 \end{gather*}
 satisfying the usual unit and associativity conditions. Take $\eps=id$ and note that $\phi$ needs to assign an isomorphism
 \begin{gather*}
 \phi_{x,y}\colon \ \LE_{\exp(x),\exp(y)} \longrightarrow \CC
 \end{gather*}
to each pair of objects $(x,y)$ of $\mT$. By construction of $\LE$, we have a trivialization of~$\LE$ over~$\mathfrak t\times\mathfrak t$, and this trivialization serves as our~$\phi$. One checks that $(F,\phi,\eps)$ satisf\/ies the axioms of a~monoidal equivalence from $\(\mT,\bullet,0\)$ to $\(\mI,\circledast,(1,\CC)\)$.
\end{proof}
\begin{Remark}\label{rem:homom}
 Construction \ref{cons:mbg} takes the sum of two bilinear forms to the tensor product of the corresponding multiplicative bundle gerbes.
\end{Remark}
\begin{Definition}
 We will write $J^t$ for the bilinear form
 \begin{gather*}
 J^t(m,n) = J(n,m).
 \end{gather*}
 We say that $J$ is {\em symmetric} if $J=J^t$ and that~$J$ is {\em skew symmetric} if $J=-J^t$. A symmetric bilinear form~$I$ is called {\em even} if
 \begin{gather}\label{eq:2Z}
 I(m,m) \in 2 \ZZ \qquad \text{for }m\in \Lat.
 \end{gather}
\end{Definition}

\begin{Proposition}\label{prop:J^t} The Lie $2$-groups $(\mT,\bullet_J,0)$ and $(\mT,\bullet_{J^t},0)$ are equivalent as Lie $2$-group extensions of $T$ by $\pt\mmod U(1)$.
\end{Proposition}
\begin{proof}
 A monoidal equivalence is given by the functor
 \begin{align*}
 F \colon \ \mT & \longrightarrow \mT, \\
 \mathfrak t & \stackrel\id\longrightarrow\mathfrak t \hspace{44.5mm} \text{on objects,}\\
 (x,m,z) & \longmapsto (x,m,z\cdot\exp(J(x,m))) \qquad\text{on arrows,}
 \end{align*}
 together with the natural transformation
 \begin{gather*}
 \phi \colon \ F(-)\bullet_{J^t} F(-) \Longrightarrow F(-\bullet_{J}-)
 \end{gather*}
 with
 \begin{gather*}
 \phi_{x,y}\colon \ x+y \xrightarrow{\,\,\,\,\,\,\exp(J(x,y))\,\,\,\,\,\,} x+y.\tag*{\qed}
 \end{gather*}\renewcommand{\qed}{}
\end{proof}

\begin{Corollary}\label{cor:skew}\quad
 \begin{enumerate}\itemsep=0pt
 \item[$(i)$] If $I$ is an even symmetric bilinear form on {\rm $\Lat$}, then the multiplicative bundle gerbe associated to $I$ possesses a square root.
 \item[$(ii)$] If $B$ is a skew symmetric integral bilinear form on~{\rm $\Lat$}, then $B$ yields a trivial $2$-group extension of $T$.
 \end{enumerate}
\end{Corollary}
\begin{proof}
 Claim (i) follows from Proposition \ref{prop:J^t}, using Remark \ref{rem:homom} and the fact
 that every even symmetric bilinear form $I$ can be
 written in the form $I=J+J^t$ for an integer-valued bilinear form
 $J$. For instance, f\/ix a basis $(b_1,\dots,b_r)$ of $\Lat$ and set
 \begin{eqnarray*}
 J(b_i,b_j) & = &
 \begin{cases}
 \frac12 I(b_i,b_i) & \text{ if }i=j,\\
 I(b_i,b_j) & \text{ if }i<j,\\
 0 & \text{ else.}
 \end{cases}
 \end{eqnarray*}
 Claim (ii) follows, similarly, from the fact that every skew symmetric bilinear form
 $B$ can be written in the form $B=J-J^t$ for an integer-valued
 bilinear form $J$.
\end{proof}
\begin{Corollary}
 Let $J$ be an integer-valued bilinear form on {\rm $\Lat$}. Then, up to equivalence, the categorical extension $(\mT,\bullet_J)$ of $T$ by $\pt\mmod U(1)$ only depends on the even bilinear form
 \begin{gather*}
 I(m,n) = J(m,n) + J(n,m).
 \end{gather*}
\end{Corollary}
\begin{proof} Let $J_1$ and $J_2$ be two integer-valued bilinear forms on $\Lat$, and assume that
 \begin{gather*}
 J_1 + J_1^t = J_2 + J_2^t.
 \end{gather*}
 Then $J_1-J_2$ is skew symmetric. By Corollary \ref{cor:skew}, it follows that the multiplicative bundle gerbe obtained from $J_1-J_2$ is
 trivial. Using Remark \ref{rem:homom}, we conclude that the multiplicative bundle gerbes obtained from $J_1$ and $J_2$ are isomorphic.
\end{proof}

\section{The example of the circle}\label{sec:circle}
Let $\Lat=\ZZ$ and $J(m,n) = mn$.
The {\em basic circle extension} $\mU(1)$ of the circle group $U(1)$ consists of the following data:
\begin{enumerate}\itemsep=0pt
 \item[(i)] the trivial bundle gerbe over $U(1)$,
 \item[(ii)] the line bundle $\LE$ on
$ U(1)\times U(1) = \RR^2/\ZZ^2$ def\/ined by the multipliers
 \begin{align*}
 f_{(m,n)}\colon \ \RR^2 & \longrightarrow U(1),\\
 (x,y) & \longmapsto \exp(my)
 \end{align*}
 for $(m,n)\in\ZZ^2$,
 \item[(iii)] the canonical isomorphism
 \begin{gather*}
 \alpha\colon \ \on{m}_{12}^*\LE\tensor {\rm pr}_{12}^*\LE \cong \on{m}_{23}^*\LE\tensor {\rm pr}_{23}^*\LE
 \end{gather*}
 over $U(1)\times U(1)\times U(1)$.
\end{enumerate}
For $k\in\ZZ$, the {\em $k^{\rm th}$ circle extension $\mU(1)_k$} of $U(1)$ is obtained by replacing the multipliers with $\exp({kmy})$.
\begin{Remark}\label{rem:Chern--Weil}
 Recall that gerbes on $U(1)$ are classif\/ied, up to stable isomorphism, by their Dixmier--Douady class in
$H^3(U(1);\ZZ) = 0$.
So, any bundle gerbe over $U(1)$ is trivializable. A multiplicative structure on the trivial bundle gerbe over $U(1)$ consists of a line bundle on $U(1)\times U(1)$ plus extra data, encoding associativity. Line bundles on $U(1)\times U(1)$ are classif\/ied, up to isomorphism, by their f\/irst Chern class in
 \begin{gather*}
 H^2(U(1)\times U(1);\ZZ) \cong \on{Alt^2}\big(\ZZ^2\big),
 \end{gather*}
 the group of skew symmetric bilinear forms on~$\ZZ^2$. This group is inf\/inite cyclic, generated by the determinant, and we claim
$ c_1(\LE) = -\det$.
 This fact is proved using Chern--Weil theory, as illustrated in Fig.~\ref{fig1}:
 \begin{figure}[h!]\centering
 \begin{tikzpicture}
 \matrix (m) [matrix of math nodes, row sep=3em, column sep=3em]
 { {\Omega^2_{U(1)\times U(1)}} & {-{\rm d}x\wedge {\rm d}y} & & \\
 {\Omega^1_{U(1)\times U(1)}} & -x{\rm d}y & -m{\rm d}y &
 \\
 {\ul\TT_{U(1)\times U(1)}} & {} & \exp(-my) &
 1 \\
 {U(1)\times U(1)} & \RR^2 & \RR^2\times\ZZ^2 &
 \RR^2\times\ZZ^2\times\ZZ^2 \\ };
 \draw[|->] (m-2-2) -- (m-1-2);
 \draw[|->] (m-3-3) -- (m-2-3);
 \draw[|->] (m-2-2) -- (m-2-3);
 \draw[|->] (m-3-3) -- (m-3-4);
 { [start chain] \chainin (m-3-1);
 \chainin (m-2-1)
 [join={node[left,labeled] {d\log}}];
 \chainin (m-1-1)
 [join={node[left,labeled] {d}}];
 }
 \draw[->, very thick] (-3.2,-3.6) -- (-3.2,3.6) ;
 \draw[dashed, very thick, red] (-6,2) -- (6,2) ;
 \draw[->, very thick] (-6,-2) -- (6,-2) ;
 \node at (0,3.8) {};
 \end{tikzpicture}
 \caption{This 2-cocycle in the truncated \v{C}ech--de Rham double complex relates the multipliers of $\LE$ to the 2-form $-\det$.} \label{fig1}
 \end{figure}
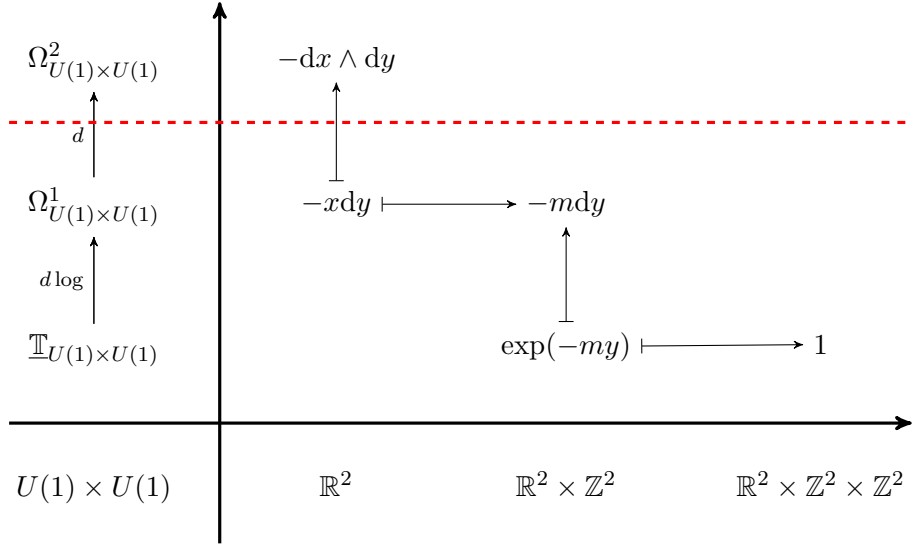

The entry $\exp(-my)$ is the transition function of the bundle $\LE$ for the cover $\RR^2\to\RR^2/\ZZ^2$. This convention, that the transition functions in \v{C}ech cohomology are taken to be the inverse multipliers, comes from algebraic geometry: with the identif\/ications $\RR^2\cong \CC$ and $\ZZ^2\cong \langle 1,\tau\rangle$, the line bundle $\LE$ is the topological line bundle underlying any degree $-1$ line bundle on the elliptic curve $\CC/\langle1,\tau\rangle$. Note also that this argument equips $\LE$ with a~connection, namely the one with 1-form
\begin{gather*}
 \omega = -x{\rm d}y
\end{gather*}
on ${\RR\times\RR}$ and with curvature
\begin{gather*}
{\rm d}\omega = -{\rm d}x\wedge {\rm d}y.
\end{gather*}
Similarly, the bundle $\LE^{\tensor k}$ with multipliers $\exp({kmy})$ has connection $kxdy$ and Chern class $k \cdot \det$.
\end{Remark}
\begin{Remark}
We refer the reader to \cite[Theorem~2.2.11]{Brylinski:Loop_spaces_characteristic_classes} for a detailed discussion of line bundles with connection in this setting. This example is not new. A construction of the multiplicative gerbe~$\mU(1)$ (but with a~dif\/ferent connection) already turns up in~ \cite[Section~3]{Brylinski:McLaughlin:The_geometry_of_degree-4_characteristic_classes_II}. A construction of $\mU(1)$ as equivariant bundle gerbe is given in \cite{Murray:Roberts:Stevenson:Vozzo}
\end{Remark}
\section{The classif\/ication}\label{sec:classification}
Recall from \cite{Schommer-Pries:Central_extensions} that, up to equivalence, the 2-group extensions of $T$ by $\pt\mmod U(1)$ are classif\/ied by degree three Lie group cohomology classes of $T$ with coef\/f\/icients in $U(1)$. There are a number of constructions of Lie group cohomology, and a
unifying axiomatic framework was recently given in~\cite{Wagemann:Wockel:A_cocycle_model}. We choose to work with the \v{C}ech-simplicial double complex of \cite{Brylinski:Loop_spaces_characteristic_classes} and \cite[Proposition~5.2]{Carey:Johnson:Murray:Stevenson:Wang}. So, the relevant
cohomology group is
\begin{gather*}
 H^3_{gp}(BT;U(1)) = \vH^3(BT_\bullet;U(1)).
\end{gather*}
The goal of this section is to analyze the composite of isomorphisms
\begin{gather}\label{eq:H^4}
 S^*(\Lambda) \cong H^{2*}(BT;\ZZ) \cong \vH^{2*}(BT_\bullet;\ZZ) \cong \vH^{2*-1}(BT_\bullet;U(1))
\end{gather}
in degree $*=2$. Here $S^*(\Lambda)$ is the symmetric algebra of the weight lattice
$ \Lambda = \Hom(\vL,\ZZ)$.
We may think of elements of
$ S^2(\Lambda) = (\Lambda\tensor\Lambda) /S_2$
as homogeneous polynomials of degree 2 in the weights, and we have the symmetrization map
\begin{align*}
 S^2(\Lambda) & \longrightarrow \on{Bil}(\vL,\ZZ),\\
 \lambda\mu & \longmapsto \lambda\tensor\mu + \mu\tensor\lambda,
\end{align*}
identifying $S^2(\Lambda)$ with the group of even symmetric bilinear forms on $\Lat$.
\begin{Theorem}\label{thm:classification} Let $I$ be an even symmetric bilinear form on {\rm $\Lat$}, and let $J$ be an integral bilinear form on {\rm $\Lat$} satisfying $I=J+J^t$. Then $I$ classifies the multiplicative bundle gerbe obtained from {\rm $(\Lat,-J)$} via Construction~{\rm \ref{cons:mbg}}.
\end{Theorem}

\begin{proof} Let $ET$ be a contractible space on which $T$ acts freely, and let $BT=ET/T$ be our model for the clasifying space of $T$. On $BT$, we
 have the line bundles
 \begin{gather*}
 L_\lambda = ET\times_T\CC_\lambda \cong \(ET\times\CC_{-\lambda}\)/T,
 \end{gather*}
 $\lambda\in\Lambda$, where $\CC_\lambda$ denotes the irreducible complex representation of $T$ with weight $\lambda$. The f\/irst isomorphism in~\eqref{eq:H^4} is
 \begin{align}
 \label{eq:Borel} S^*(\Lambda) & \stackrel\cong\longrightarrow H^{2*}(BT;\ZZ),\\
 \notag \lambda & \longmapsto c_1(L_\lambda).
 \end{align}
 Its def\/inition goes back to Borel, and it is an isomorphism of graded rings. It is, of course, well known that the odd cohomology groups vanish, but this fact will not concern us here. To def\/ine the second isomorphism in \eqref{eq:H^4}, we use the homeomorphism
 \begin{align*}
 ET \times T \times \dots\times T & \stackrel\cong\longrightarrow ET\times_{BT}\dots\times_{BT} ET,\\
 (e,t_1,\dots,t_n) & \longmapsto (e,t_1e,\dots,t_ne)
 \end{align*}
 to interpret the maps
 \begin{align*}
 ET & \longrightarrow BT, \\
 ET \times \mathfrak t & \longrightarrow ET\times_{BT} ET
 \end{align*}
 as a hypercover of $BT$. The \v{C}ech double complex of this hypercover can be identif\/ied with the \v{C}ech simplicial double complex $\vC^*(BT_\bullet;\ZZ)$. Under this identif\/ication, the cup product becomes
 \begin{gather*}
 (f\cup g) (x_0,\dots,x_{r+s}) = ({\rm pr}_{1}^*f)(x_0,\dots x_r)\cdot ({\rm pr}_{2}^*g)(x_r,\dots x_{r+s}),
 \end{gather*}
 where $r$ and $s$ are the \v{C}ech degrees of $f$ and $g$, while ${\rm pr}_1$ is the projection onto the f\/irst $\deg_{\rm simp}(f)$ factors, and
 ${\rm pr}_2$ is the projection onto the last $\deg_{\rm simp}(g)$ factors. To determine the image of~$\lambda\mu$ in $H^4(BT_\bullet;\ZZ)$,
 we determine the images of $\lambda$ and $\mu$ and then take the cup product. The f\/irst Chern class of~$L_\lambda$ in \v{C}ech hypercohomology is
 obtained by applying $\delta_{\text{\v{C}ech}}\log$ to the inverse multipliers
 \begin{align*}
 ET\times \mathfrak t & \longrightarrow U(1),\\
 (\eta,x) & \longmapsto \exp(\lambda(x))
 \end{align*}
 of $L_\lambda$.
 So, $c_1(L_\lambda)$ is represented by the degree $(1,1)$ cocycle
 \begin{align*}
 \mathfrak t\times\vL & \longrightarrow \ZZ, \\
 (x,m) & \longmapsto \lambda(m)
 \end{align*}
 in the \v{C}ech-simplicial double complex. Given two weights, $\lambda$ and $\mu$, their cup product is represented by the cocycle
 \begin{align*}
 \mathfrak t^2\times(\vL)^2\times(\vL)^2 & \longrightarrow \ZZ, \\
 \left(\left( \begin{matrix}
 x\\y
 \end{matrix}\right),
 \left( \begin{matrix}
 m\\n
 \end{matrix}\right),
 \left(\begin{matrix}
 k\\l
 \end{matrix}\right)\right)
 & \longmapsto \lambda(k)\cdot\mu(n)
 \end{align*}
 in $\vC^2(BT_2;\ZZ)$. The last isomorphism in \eqref{eq:H^4} is the inverse of the connecting homomorphism for the short exact sequence of presheaves
 \begin{gather*}
 0 \longrightarrow \ZZ \longrightarrow \ul\RR \longrightarrow \ul{U(1)} \longrightarrow 0.
 \end{gather*}
 Let $J=\mu\tensor\lambda$. As in Remark~\ref{rem:Chern--Weil}, the transition function for $\LE_{-J}$ is the inverse multiplier. So, Construction~\ref{cons:mbg} associates to $-J$ the multiplicative bundle gerbe corresponding to the $U(1)$-valued \v{C}ech-simplicial
 3-cocycle $(1,\lambda(m)\mu(y),1,1)$. The image of this cocycle under the connecting homomorphism is the integral
 4-cocycle \begin{gather*}(1,1,\lambda(k)\mu(n),1,1)=c_1\(L_\lambda\)\cdot c_1\(L_\mu\).\end{gather*}

 \begin{figure}[t]\centering
 \begin{tikzpicture}[scale=.8]
 \matrix (m) [matrix of math nodes, row sep=2.4em, column sep=3em]
 { {T\times T\times T} & & &
 \\
 {T\times T} & & {\lambda(m)\mu(y)} & \lambda(k)\mu(n)
 \\
 {T} & {\lambda(x)} &\lambda(m) &
 \\
 {\pt} & {} & &
 \\
 {} & \mathfrak t^\bullet & (\mathfrak t\times\vL)^\bullet &
 (\mathfrak t \times\vL\times\vL)^\bullet \\ };
 \draw[|->] (m-2-3) -- (m-2-4) node [midway, above] {$\delta$};
 \draw[|->] (m-3-2) -- (m-3-3) node [midway, above] {$\delta$};
 \draw[->, very thick] (-3,-4) -- (-3,4.5) ;
 \draw[->, very thick] (-6,-3) -- (7.5,-3) ;
 \end{tikzpicture}
\caption{$\delta_{\text{\v{C}ech}}\colon \check{C}^1(BT_2;\ul\RR) \longrightarrow \check{C}^2(T_2;\ul\RR)$.}\label{fig2}
 \end{figure}

 The general case now follows from Remark~\ref{rem:homom} together with the fact that any choice of $J$ can be written as linear combination of summands of the form $\mu\tensor\lambda$.
\end{proof}
\begin{Remark} In particular, we have seen that the underlying gerbe of any multiplicative bundle gerbe on $T$ is trivial. This is consistent with work of Waldorf
 \cite[Proposition~2.10]{Waldorf:Multiplicative_bundle_gerbes_with_connection}, who
 identif\/ies the forgetful map
 \begin{gather*}
 \{\text{multiplicative bundle gerbes on }G\}/{\cong} \longrightarrow \{\text{bundle gerbes on }G\}/{\cong}
 \end{gather*}
 with the inverse transgression
$\tau\colon H^4(BG;\ZZ) \longrightarrow H^3(G;\ZZ)$
 in the Leray--Serre spectral sequence for the f\/ibration $EG\to BG$.
 For a compact, connected Lie group $G$, this inverse transgression map was calculated (in all degrees) by Chern and Simons. In the relevant degrees, their result is summarized by the commuting diagram\footnote{The formula
 \cite[equation~(3.10)]{Chern:Simons:Characteristic_forms_and_geometric_invariants}
 is often cited in its original form
 \begin{gather*}
 I \longmapsto -\tfrac1{12} I(\omega\wedge[\omega,\omega]),
 \end{gather*}
 where $\omega$ is the right-invariant Maurer--Cartan form on $G$. A look at the def\/initions on p.~50 of~\cite{Chern:Simons:Characteristic_forms_and_geometric_invariants} identif\/ies $I(\omega\wedge[\omega,\omega])$ with the bi-invariant
 3-form on $G$ whose restriction to $\mathfrak g=T_1G$ equals~$12 \nu$.}
 \begin{equation*}
 \begin{tikzpicture}
 \node at (0,3) [name=A] {$H^4(BG;\RR)$};
 \node at (4,3) [name=B] {$H^3(G;\RR)$};
 \node at (0,0) [name=C] {$\(S^2\mathfrak g^*\)^{\rm Ad}$};
 \node at (4,0) [name=D] {$\(\Lambda^3\mathfrak g^*\)^{\rm Ad}$};
 \node at (0,-1) [name=E,text width=1cm] {\quad$I$};
 \node at (4,-1) [name=F, text width=1.5cm] {\quad$-\nu$,};
 \draw[->] (A) -- (B) node[midway,above]{$\tau$};
 \draw[->] (C) -- (A) node[midway,right]{$\cong$}
 node[midway,left, text width= 2.5cm]{Chern--Weil isomorphism};
 \draw[->] (B) -- (D) node[midway,left]{$\cong$}
 node[midway,right, text width= 3.8cm]{\quad Chevalley--Eilenberg\\
 \quad isomorphism};
 \draw[->] (C) -- (D);
 \draw[|->] (E) -- (F);
 \end{tikzpicture}
 \end{equation*}
 where $\nu$ is the Cartan 3-form
 associated to $I$,
 \begin{gather*}
 \nu(x,y,z) = I([x,y],z).
 \end{gather*}
 So, the multiplicative bundle gerbe on $G$ classif\/ied by $I$ has as its Dixmier--Douady class the left-invariant 3-form on $G$ whose restriction to~$\mathfrak g$ is~$-\nu$. Since the Lie bracket on a torus is zero, it follows that the Dixmier--Douady class of any multiplicative bundle gerbe on~$T$ vanishes.
\end{Remark}

\subsection{Connections}
Let $\LE$ be the line bundle of Construction~\ref{cons:mbg}. If $J=\lambda\tensor\mu$, then we have the connection $\nabla$ on~$\LE$ whose 1-form on ${\mathfrak t\times\mathfrak t}$ is
\begin{gather*}
 \omega = -\lambda {\rm d}\mu.
\end{gather*}
It has curvature 2-form
\begin{gather*}
 \kappa = -{\rm d}\lambda\wedge {\rm d}\mu.
\end{gather*}
For arbitrary $J$, we introduce the maps
\begin{align*}
 J^\sharp(x)\colon \ \mathfrak t & \to \RR, \\
 v & \longmapsto J(x,v).
\end{align*}
Then $\nabla$ is def\/ined by the 1-form
\begin{gather*}
 \omega_{(x,y)} = -{\rm d}(J^\sharp(x))_y.
\end{gather*}
The pair $(\LE,\nabla)$ turns the trivial bundle gerbe $\mI$ on $T$ into a multiplicative bundle gerbe with connection (in the sense of~\cite{Waldorf:Multiplicative_bundle_gerbes_with_connection}, with the remaining data trivial).

\section{Examples}\label{sec:examples}
\subsection{Maximal tori}\label{sec:maximal_tori}
Let $G$ be a simple and simply connected compact connected Lie group with maximal torus $T$ and Weyl group $W$. Then we have
\begin{gather*}
 \ZZ \cong H^3_{\rm gp}(G;U(1)) \cong H^4(BT;\ZZ)^W \cong\big( \operatorname{Bil}_{\rm ev}(\Lat;\ZZ)^{S_2}\big)^W.
\end{gather*}
The elements of this group are multiples of the Killing form, and the positive def\/inite gene\-ra\-tor~$I_{\rm basic}$ classif\/ies the Lie 2-group extension of $G$ denoted $\on{String}(G)$ in \cite{Schommer-Pries:Central_extensions}. We arrive at a~recognition principle for these 2-groups: the extension $\mG$ is equivalent to $\on{String}(G)$ if and only if its restriction to $T$ is equivalent to the categorical torus classif\/ied by $(\vL,I_{\rm basic})$.

\subsection{The Leech lattice}\label{sec:Leech}
 Another interesting example is given by the Leech lattice $\Lat=\Lambda_{\rm Leech}$ inside $\RR^{24}$, together with the standard symmetric bilinear form $I$. The group of linear isometries $O(\Lat,I)$ of the Leech lattice is the {\em Conway group} ${\rm Co}_0$. In analogy to the previous example, we view~$I$ as a ${\rm Co}_0$-invariant cohomology element,
 \begin{gather*}
 I\in H^4(BT,\ZZ)^{{\rm Co}_0},
 \end{gather*}
 where $T=\RR^{24}/\Lambda_{\rm Leech}$ is the Leech torus. We arrive at a categorical extension $\mT_{\rm Leech}$ of the Leech torus, on which ${\rm Co}_0$ acts by autoequivalences. It is now understood \cite{Johnson-Freyd:Treumann} that the Conway group has a universal categorical central extension by the cyclic group with $24$ elements, which we believe to be closely related, but not equal to the symmetries of the categorical Leech torus.
\subsection{Niemeyer lattices}
 Similarly, we can choose $\Lat$ as one of the Niemeyer lattices, and $I$ as the standard symmetric bilinear form on $\RR^{24}$, making $\Lat$ an even unimodular lattice. If $\Lat$ equals $A_1^{24}$ or $A_2^{12}$, then the group of linear isometries of $(\Lat,I)$ is the {\em Mathieu group} $M_{24}$ (respectively $M_{12}$), and we have
 \begin{gather*}
 I\in H^4(BT,\ZZ)^{M_{24}} \qquad\text{respectively}\qquad I\in H^4(BT,\ZZ)^{M_{12}},
 \end{gather*}
 classifying two categorical tori on which the respective Mathieu groups act by autoequivalences.

\section{Extraspecial categorical 2-groups}\label{sec:Extraspecial}
 Assume we are given a strict action of a group $G$ by functors on a category $\mC$. In other words, assume that $\mC$ comes equipped with endofunctors $\varrho(g)$, one for each $g\in G$, satisfying
 \begin{gather*}
 \varrho(g)\varrho(h)=\varrho(gh)\qquad \text{and}\qquad \varrho(1)={\rm id}.
 \end{gather*}
 In this situation, we will allow ourselves to drop $\varrho$ from the notation and simply write $g$ for the functor $\varrho(g)$.
\begin{Definition}[Grothendieck]
 An {\em equivariant object} of $\mC$ is a pair
 \begin{gather*}(x,\mathbf e)=(x,\{e_g\}_{g\in G}),\end{gather*} where $x\in {\rm ob}(\mC)$ and the $\isomap{e_g}xgx$ are isomorphisms in $\mC$ satisfying
 \begin{gather} \label{eq:equivariant_object}
 (ge_h)\circ e_g = e_{gh},
 \end{gather}
 for all $g$ and $h$ in $G$. An {\em equivariant arrow} from $(x,\mathbf e)$ to $(y,\mathbf f)$ is an arrow $\map axy$ in $\mC$ that is compatible with $\mathbf e$ and $\mathbf f$ in the following sense
 \begin{gather*}
 f_g\circ a = (ga)\circ e_g,
 \end{gather*}
 for all $g\in G$.
\end{Definition}
\begin{Definition} We will refer to the category $\mC^G$ of equivariant objects of $\mC$ and equivariant arrows between them as the {\em categorical fixed points} $\mC^G$ of the action of $G$ on $\mC$.
\end{Definition}

The goal of this section is to identify the f\/ixed point category of the action of $\{\pm1\}$ on $\mT$, where $-1$ acts by the auto-equivalence associated to the crossed module automorphism
\begin{equation*}
 \begin{tikzpicture}
 \node at (-1.5,-2) {$\operatorname{inv}\colon$};
 \node at (0,0) [name=a] {$(m,z)$};
 \node at (3,0) [name=b] {$(-m,z)$};
 \node at (0,-1) [name=c] {$\Lat\times U(1)$};
 \node at (3,-1) [name=d] {$\Lat\times U(1)$};
 \node at (0,-3) [name=e] {$\mathfrak t$};
 \node at (3,-3) [name=f] {$\mathfrak t$};
 \node at (0,-4) [name=g] {$x$};
 \node at (3,-4) [name=h] {$-x$.};
 \draw[|->] (a) -- (b);
 \draw[->] (c) -- (d);
 \draw[->] (e) -- (f);
 \draw[->] (c) -- (e);
 \draw[->] (d) -- (f);
 \draw[|->] (g) -- (h);
 \end{tikzpicture}
\end{equation*}
Explicitly, the functor $\varrho(-1)$ sends the object $x$
to $-x$ and the arrow \makebox{$x\xrightarrow{\,\,\,\,z\,\,\,\,}x+m$}
to the arrow
\makebox{$-x\xrightarrow{\,\,\,\,z\,\,\,\,}-x-m$}.
Since this is an action by strictly monoidal functors, the categorical f\/ixed points $\mT^{\{\pm1\}}$ inherit the structure of a strict categorical group. Consider now the points of order two in $T$. These form the elementary abelian 2-group\footnote{In the context of this section, the term 2-group is used in the sense of $p$-group, $p$ a prime, not in the
 sense of categorical group.}
\begin{gather*}
 T[2] \,\cong \,\Lat\tensor\,\FF_2.
\end{gather*}
We have a central extension $\widetilde{T[2]}$ of $T[2]$,
def\/ined by the $\FF_2$-valued 2-cocycle
\begin{gather*}
 J_{\FF_2}\,:=\,J\tensor \,\FF_2.
\end{gather*}
Central extensions of this form are know as {\em extraspecial $2$-groups} and classif\/ied by their {\em Arf invariant}, i.e., by
the quadratic form
\begin{gather*}
 v \longmapsto J_{\FF_2}(v,v)
\end{gather*}
on the $\FF_2$-vector space $\Lat\!\tensor\FF_2$. Writing its center multiplicatively,
our central extension takes the form
\begin{equation*}
\begin{tikzpicture}
 \node at (0.2,0) [name=a] {$1$};
 \node at (2,0) [name=b] {$\{\pm1\}$};
 \node at (4.3,0) [name=c] {$\widetilde{T[2]}$};
 \node at (6.5,0) [name=d] {$T[2]$};
 \node at (8.4,0) [name=e] {$1$.};
 \draw[->] (a) -- (b);
 \draw[->] (b) -- (c);
 \draw[->] (c) -- (d);
 \draw[->] (d) -- node [midway, above] {} (e);
\end{tikzpicture}
\end{equation*}
Similarly, we have the (non-canonically trivial) central extension
\begin{equation*}
\begin{tikzpicture}
 \node at (0.2,0) [name=a] {$1$};
 \node at (2,0) [name=b] {$\{\pm1\}$};
 \node at (4.3,0) [name=c] {$\widetilde{\Lat}$};
 \node at (6.5,0) [name=d] {$\Lat$};
 \node at (8.4,0) [name=e] {$1$};
 \draw[->] (a) -- (b);
 \draw[->] (b) -- (c);
 \draw[->] (c) -- (d);
 \draw[->] (d) -- node [midway, above] {} (e);
\end{tikzpicture}
\end{equation*}
with 2-cocycle
\begin{gather*}
 (m,n) \longmapsto (-1)^{J(m,n)}.
\end{gather*}
\begin{Theorem}
 The crossed module corresponding to the strict categorical group $\mT^{\{\pm1\}}$ is isomorphic to the crossed module
 \begin{align*}
 \Sigma\colon \ \text{\rm $\Lat$}\times U(1) & \longrightarrow \widetilde{\text{\rm $\Lat$}} ,\\
 (m,z) & \longmapsto (2m,1),
 \end{align*}
 where $(n,\eps)\in \widetilde{\text{\rm $\Lat$}} $ acts on $\text{\rm $\Lat$}\times U(1)$ by
 \begin{gather*}(m,z)\longmapsto \big(m,z\cdot(-1)^{J(m,n)}\big).\end{gather*}
 In particular, $\mT^{\{\pm1\}}$ is part of a categorical central extension
 \begin{equation*}
 \begin{tikzpicture}
 \node at (-0.2,0) [name=a] {$1$};
 \node at (1.8,0) [name=b] {$\pt\mmod U(1)$};
 \node at (4.4,0) [name=c] {$\mT^{\{\pm1\}}$};
 \node at (6.5,0) [name=d] {$\widetilde{T[2]}$};
 \node at (8.2,0) [name=e] {$1$.};
 \draw[->] (a) -- (b);
 \draw[->] (b) -- (c);
 \draw[->] (c) -- (d);
 \draw[->] (d) -- node [midway, above] {} (e);
\end{tikzpicture}
\end{equation*}
\end{Theorem}
\begin{proof} For $\mT^{\{\pm1\}}$, condition \eqref{eq:equivariant_object} reads \begin{gather*}e_1=\id_x \qquad \text{and} \qquad {\rm inv} (e_{-1})\circ e_{-1} = \id_x.\end{gather*}
 So, we have a fully faithful and strictly monoidal embedding
 \begin{equation*}
 \begin{tikzpicture}
 \node at (0,0) [name=A, anchor=east] {$\mT^{\{\pm1\}}$};
 \node at (1,0) [name=B, anchor=west] {$\mT^I$,};
 \node at (0,-0.8) [name=C, anchor=east] {$(x,\boldsymbol e)$};
 \node at (1,-0.8) [name=D, anchor=west] {$e_{-1}$,};
 \draw[->] (A) -- (B);
 \draw[|->] (C) -- (D);
 \end{tikzpicture}
 \end{equation*}
 where $\mT^I$ denotes the arrow category of $\mT$. This embedding identif\/ies the objects of $\mT^{\{\pm1\}}$ with the image of the
 injective group homomorphism
 \begin{equation*}
 \begin{tikzpicture}
 \node at (0,0) [name=A, anchor=east] {$\widetilde\Lat$};
 \node at (1,0) [name=B, anchor=west] {${\rm ob}\big(\mT^I\big)$,};
 \node at (0,-0.8) [name=C, anchor=east] {$(m,\eps)$};
 \node at (1,-0.8) [name=D, anchor=west]
 {$\big({-}\frac
 m2\xrightarrow{\,\,\,\,\eps\,\,\,\,}\frac m2\big)$.};
 \draw[->] (A) -- (B);
 \draw[|->] (C) -- (D);
 \end{tikzpicture}
 \end{equation*}
 The kernel of the source of $\mT^{\{\pm1\}}$ is identif\/ied with the image
 of the group homomorphism
 \begin{equation*}
 \begin{tikzpicture}
 \node at (0,0) [name=A, anchor=east] {$\Lat\times U(1)$};
 \node at (1,0) [name=B, anchor=west] {${\rm arrows}\big(\mT^I\big)$,};
 \node at (0,-1.5) [name=C, anchor=east] {$(m,\eps)$};
 \node at (1,-1.5) [name=D, anchor=west]
 {\begin{tikzpicture}
 \node (X) {$0$};
 \node [below=of X] (Y) {$0$};
 \node [right=of X] (Z) {$-m$};
 \node [below=of Z] (W) {$m$,};
 \draw[double equal sign distance] (X) -- (Y);
 \draw[->] (X) -- node [midway, above] {$z$} (Z);
 \draw[->] (Y) -- node [midway, above] {$z$} (W);
 \draw[->] (Z) -- node [midway, right] {$1$} (W);
 \end{tikzpicture}
 };
 \draw[->] (A) -- (B);
 \draw[|->] (C) -- (D);
 \end{tikzpicture}
 \end{equation*}
 and under these identif\/ications, the target map becomes the homomorphism $\Sigma$ of the theorem. It remains to identify the conjugation action of the objects on the
 kernel of the source.
We have
\begin{equation*}
 \begin{tikzpicture}[scale=1.1]
 \node at (0.6,2) [name=a] {$0$};
 \node at (2.6,2) [name=g] {$-m$};
 \node at (4,2) [name=c] {$-\frac n2$};
 \node at (6,2) [name=i] {$-\frac n2$};
 \node at (8,2) [name=e] {$-\frac n2$};
 \node at (11.5,2) [name=k] {$-\frac{2m+n}2$};
 \node at (0.6,0) [name=b] {$0$};
 \node at (2.6,0) [name=h] {$m$};
 \node at (4,0) [name=d] {$\frac n2$};
 \node at (6,0) [name=j] {$\frac n2$};
 \node at (8,0) [name=f] {$\frac n2$};
 \node at (11.5,0) [name=l] {$\frac{2m+n}2$};
 \node at (3.3,1) {$\bullet$};
 \node at (7,1) {$=$};
 \draw[double equal sign distance] (a) -- (b);
 \draw[->] (c) -- node[midway,right] {$\eps$} (d);
 \draw[->] (e) -- node[midway,right] {$\eps$} (f);
 \draw[->] (a) -- node[midway,above] {$z$} (g);
 \draw[->] (b) -- node[midway,above] {$z$} (h);
 \draw[double equal sign distance] (c) -- (i);
 \draw[double equal sign distance] (d) -- (j);
 \draw[->] (e) -- node[midway,above] {\!\small{$z(-1)^{J(m,n)}$}} (k);
 \draw[->] (f) -- node[midway,above] {\!\small{$z(-1)^{J(m,n)}$}} (l);
 \draw[->] (g) -- node[midway,right] {$1$} (h);
 \draw[->] (i) -- node[midway,right] {$\eps$} (j);
 \draw[->] (k) -- node[midway,right] {$\eps$} (l);
 \end{tikzpicture}
\end{equation*}
and
\begin{equation*}
 \begin{tikzpicture}[scale=1.1]
 \node at (4,2) [name=a] {$0$};
 \node at (6,2) [name=g] {$-m$};
 \node at (.6,2) [name=c] {$-\frac n2$};
 \node at (2.6,2) [name=i] {$-\frac n2$};
 \node at (8,2) [name=e] {$-\frac n2$};
 \node at (10.5,2) [name=k] {$-\frac{2m+n}2$};
 \node at (4,0) [name=b] {$0$};
 \node at (6,0) [name=h] {$m$};
 \node at (.6,0) [name=d] {$\frac n2$};
 \node at (2.6,0) [name=j] {$\frac n2$};
 \node at (8,0) [name=f] {$\frac n2$};
 \node at (10.5,0) [name=l] {$\frac{2m+n}2$.};
 \node at (3.3,1) {$\bullet$};
 \node at (7,1) {$=$};
 \draw[double equal sign distance] (a) -- (b);
 \draw[->] (c) -- node[midway,right] {$\eps$} (d);
 \draw[->] (e) -- node[midway,right] {$\eps$} (f);
 \draw[->] (a) -- node[midway,above] {$z$} (g);
 \draw[->] (b) -- node[midway,above] {$z$} (h);
 \draw[double equal sign distance] (c) -- (i);
 \draw[double equal sign distance] (d) -- (j);
 \draw[->] (e) -- node[midway,above] {{$z$}} (k);
 \draw[->] (f) -- node[midway,above] {$z$} (l);
 \draw[->] (g) -- node[midway,right] {$1$} (h);
 \draw[->] (i) -- node[midway,right] {$\eps$} (j);
 \draw[->] (k) -- node[midway,right] {$\eps$} (l);
 \end{tikzpicture}
\end{equation*}
So, the conjugation action is as claimed.
\end{proof}

\begin{Remark} Let $I=J+J^t$, and consider the integer-valued quadratic form
 \begin{gather*}
 \phi(m) = \tfrac12 I(m,m)
 \end{gather*}
on $\Lat$. Its reduction mod 2 is the Arf invariant of $\widetilde{T[2]}$. The form $I$ can be recovered from $\phi$ by the identity
 \begin{gather*}
 I = \delta_{\rm simp}(\phi),
 \end{gather*}
i.e.,
 \begin{gather*}
 I(m,n) = \phi(m+n)-\phi(m)-\phi(n).
 \end{gather*}
\end{Remark}

\begin{Example}
There is a prominent
subgroup $C$ of the Monster, sitting in a non-split extension
\begin{gather*}
\begin{tikzpicture}
 \node at (0.2,0) [name=a] {$1$};
 \node at (2,0) [name=b] {$\widetilde{T[2]}$};
 \node at (4.3,0) [name=c] {$C$};
 \node at (6.5,0) [name=d] {${\rm Co}_1$};
 \node at (8.4,0) [name=e] {$1$,};
 \draw[->] (a) -- (b);
 \draw[->] (b) -- (c);
 \draw[->] (c) -- (d);
 \draw[->] (d) -- node [midway, above] {} (e);
\end{tikzpicture}
\end{gather*}
where ${\rm Co}_1$ is the Conway group ${\rm Co}_1={\rm Co}_0/{\pm1}$. This subgroup is typically the f\/irst step in the construction of the Monster, see for instance~\cite{Tits:Le_Monstre} or~\cite{Conway:Sloane:Sphere_packings}.
\end{Example}

\section{From categorical groups to loop groups}\label{sec:loops}
The relationship between loop groups and categorical Lie groups has been the subject of extensive study. In~\cite{Baez:Stevenson:Crans:Schreiber}, the authors use central extensions of loop groups to def\/ine 2-groups. We will mostly be interested in the other direction and work with the transgression-regression machine of~\cite{Waldorf:A_construction_of_string_2-group_models,Waldorf:Transgression_to_loop_spaces_and_its_inverse_II}. In this context, the term `{\em transgression}' refers to a recipe for turning multiplicative bundle gerbes with connection into central extensions of loop groups. Let $\mL T$ be the group of piecewise smooth loops in $T$. This itself is not a manifold, but it is densely isomorphic to an increasing union of manifolds of loops which are smooth over a~given subdivision of~$\bbS^1$ into intervals. Here (and similarly below), we will follow Brylinski \cite[p.~96]{Brylinski:Loop_spaces_characteristic_classes} and call a map on~$\mL T$ smooth if its restriction to each such manifold is smooth. Applied to the trivial bundle gerbe~$\mI$ on $T$ equipped with the trivial connection, transgression yields the trivial principal $U(1)$-bundle over~$\mL T$. If we further apply transgression to our multiplicative structure with connection~$((\LE,\nabla),\alpha)$, we obtain the central extension of~$\mL T$ whose 2-cocycle $c$ is given by the holonomy of~$(\LE,\nabla)$. So, if $\phi$ and $\gamma$ are loops in $T$, and $(f,g)$ is any choice of lift of $(\varphi,\gamma)$ to $\mathfrak t\times\mathfrak t$, then we have
 \begin{gather}
 c(\varphi,\gamma) = {\rm Hol}_{\(\LE,\nabla\)}((\varphi,\gamma))
 = \exp\(\(\int_0^1 J(f(t),\dot{g}(t)) {\rm d}t\) - J\(\Delta_f,g(0)\)\), \label{eq:cocycle}
 \end{gather}
where $\Delta_f=f(1)-f(0)$.
This formula results from a variation of \cite[Proposition~6.1.3]{Brylinski:Loop_spaces_characteristic_classes}. The underlying principal bundle of this central extension is trivial.
\begin{Lemma}\label{lemma7.1}
 The central extension $\widetilde{\mL T}$ defined by $c$ is isomorphic to that in {\rm \cite[Section~4.8]{Pressley:Segal}}.
\end{Lemma}
\begin{proof}
 Pressley and Segal denote our $I$ by $\langle-,-\rangle$, their $b$ can be taken to be our $J$, and they write $\Lambda$ for the
 cocharacter lattice \begin{gather*}\vT=\Hom(U(1),T).\end{gather*} Our proof follows closely
 that of \cite[Proposition~4.8.3]{Pressley:Segal}. We f\/irst note that $c$ describes the correct extension of this lattice
 \begin{gather*}
 \Lat\cong \vT \subset \mL T,
 \end{gather*}
 namely
 \begin{gather*}
 c_\Lat(m,n) = (-1)^{J(m,n)}.
 \end{gather*}
 The commutator
 $\widetilde\varphi\cdot\widetilde\gamma\cdot\widetilde\varphi\inv\cdot
 \widetilde\gamma\inv$ in $\widetilde{\mL T}$ equals
 \begin{gather*}
 \frac{c(\varphi,\gamma)}{c(\gamma,\varphi)}
 = \exp\(\int_0^1 J(f(t),\dot g(t))\, {\rm d}t -\int_0^1J(g(t),\dot f(t))\,{\rm d}t-J\(\Delta_f,g(0)\)+J(\Delta_g,f(0))\)\\
\hphantom{\frac{c(\varphi,\gamma)}{c(\gamma,\varphi)}}{}
 = \exp\(\int_0^1 I(f(t),\dot g(t))\,{\rm d}t -
 \left[J(g(t),f(t))\right]^1_0-J\(\Delta_f,g(0)\)+J(\Delta_g,f(0))\)\\
\hphantom{\frac{c(\varphi,\gamma)}{c(\gamma,\varphi)}}{} = \exp\(\int_0^1 I(f(t),\dot g(t))\,{\rm d}t - J(g(1),\Delta_f) -
 J\(\Delta_f,g(0)\) \)\\
\hphantom{\frac{c(\varphi,\gamma)}{c(\gamma,\varphi)}}{}
 = \exp\(-\int_0^1 I(\dot f(t),g(t))\,{\rm d}t +
 \left[J(f(t),g(t))\right]^1_0-J\(\Delta_f,g(0)\)+J(\Delta_g,f(0))\)\\[+4pt]
 \hphantom{\frac{c(\varphi,\gamma)}{c(\gamma,\varphi)}}{} = \exp\(-\int_0^1 I(\dot f(t),g(t))\,{\rm d}t + J\(f(1),\Delta_g\)
 + J\(\Delta_g,f(0)\)\).
 \end{gather*}
 The second and fourth equality are obtained from the f\/irst line through integration by parts. Over the identity component of $\mL T$, our extension is completely described by its Lie-algebra cocycle,
 \begin{gather*}
 \omega = \int_0^1 I(f(t),\dot g(t)) {\rm d}t
 \end{gather*}
 (from line 3 and \cite[Section~2.1]{Pressley:Segal}), which
 agrees with the expression for $\omega$ in \cite[Proposition~4.2.2]{Pressley:Segal}. The adjoint action of $\varphi\in\mL T$ on $\widetilde{L\mathfrak t}$ can be read of\/f from the last line. It sends $(g,r)\in\widetilde{L\mathfrak t}$ to
 \begin{gather*}
 \(g,r - \int_0^1I(\dot f(t),g(t)){\rm d}t\). \end{gather*}
 This agrees with the expression in \cite[Proposition~4.3.2]{Pressley:Segal}.
\end{proof}
\begin{Remark} This def\/inition of $\widetilde{\mL T}$, as transgression of a categorical torus, leads to a con\-si\-derable simplif\/ication of the picture in the standard literature on loop groups \cite[Section~4.8]{Pressley:Segal}. The expression \eqref{eq:cocycle} is less complicated than \cite[Proposition~4.8.3]{Pressley:Segal}. It is, by construction, invariant under the action of $\on{Dif\/f^+}(\bbS^1)$. There is no need to restrict ourselves to maximal tori of simply laced groups. Further, the cocycle $c$ satisf\/ies the {\em fusion rule}
 \begin{gather}\label{eq:fusion}
 c(\varphi_1*\overline\varphi_2,\gamma_1*\overline\gamma_2) \cdot
 c(\varphi_2*\overline\varphi_3,\gamma_2*\overline\gamma_3) =
 c(\varphi_1*\overline\varphi_3,\gamma_1*\overline\gamma_3).
 \end{gather}
 Here $*$ stands for concatenation of paths, and $(\varphi_1,\varphi_2,\varphi_3)$ and $(\gamma_1,\gamma_2,\gamma_3)$ are triples of paths with identical start and end points, i.e.,
 \begin{gather*}\varphi_1(0)=\varphi_2(0)=\varphi_3(0)\qquad\text{and}\qquad
 \varphi_1(1)=\varphi_2(1)=\varphi_3(1),\end{gather*} and likewise for the $\gamma_i$.
\end{Remark}

Our fourth construction of $\mT$ applies Waldorf's regression machine to reconstruct our categorical torus from the cocycle $c$.
\begin{Construction}\label{cons:loops} Given $(\Lat,J)$, we let $T$ be the torus $\Lat\tensor U(1)$ and $\widetilde{\mL T}$ the central extension of the loop group $\mL T$ def\/ined in \eqref{eq:cocycle}. We write $\mP_1T$ for the space of piecewise smooth paths based at $1$ in $T$ and
 \begin{gather*}
 \Omega T\cong \mP_1T\times_T\mP_1T
 \end{gather*}
 for the group of piecewise smooth loops based at $1$. Then we have the groupoid
 \begin{gather*}
 \begin{tikzpicture}
 \node at (1.2,0) [name=a, anchor=west] {$\widetilde{\Omega T}\big)$,};
 \node at (0,0) [name=b, anchor=east] {$\mG := \big(\mP_1T$};
 \draw[->] (1.2,.1) -- (0,.1);
 \draw[->] (1.2,-.1) -- (0,-.1);
 \end{tikzpicture}
 \end{gather*}
 where the composition of arrows is
 \begin{gather*}
 (z,\gamma_2,\gamma_3)\circ (w,\gamma_1,\gamma_2) = (zw,\gamma_1,\gamma_3).
 \end{gather*}
 The fusion rule \eqref{eq:fusion} ensures that this composition is a~group homomorphism, making $\mG$ a~grou\-poid in groups.
\end{Construction} It is probably possible to articulate the manner in which this is an inf\/inite-dimensional ind Lie-groupoid, but we will not bother, nor will we dive into the technicalities of dif\/feological spaces as in Waldorf. We are grateful to David Roberts for explaining the following way to work around smoothness issues: Fix an order preserving smooth map
\begin{gather*}
 f\colon \ [0,1] \longrightarrow [0,2]
\end{gather*}
sending $\frac12$ to $1$. We ask that all derivatives of $f$ vanish at $0$, at $\frac12$, and at $1$ and that, away from these points, $f$ is a~dif\/feomorphism onto $(0,1)\cup(1,2)$. Consider a pair of smooth paths $\varphi,\gamma \colon [0,1]\longrightarrow T$ with $\varphi(1)=\gamma(0)$. The {\em $($mollified$)$ $f$-concatentation} of $\varphi$ and $\gamma$ is def\/ined as
\begin{gather*}
 \varphi *_f \gamma = \begin{cases}
 \varphi(f(t)) & t\leq \frac12,\\
 \gamma(f(t)-1) & t\geq \frac12.
 \end{cases}
\end{gather*}
Once can now proceed as before, with $*$ replaced by $f$-concatenation, which is smoothly homotopic to ordinary concatenation but also a smooth map on the smooth function spaces. The reparametrization invariance of the cocycle $c$ ensures that the calculations remain the same. The space of composable pairs of smooth paths in Construction \ref{cons:loops} is now identif\/ied with loops that are smooth away from possibly the basepoint and $\frac12$. This is densly isomorphic to the space of smooth loops $\Omega^{\rm sm} T$. Moreover, the latter is a smooth retract of the former, as Lie groups. Waldorf's machine guarantees that Construction \ref{cons:loops} recovers our categorical torus. We give an explicit equivalence in three steps.
\begin{Lemma} The strict categorical Lie group $\mG$ of Construction~{\rm \ref{cons:loops}} corresponds to the crossed module
 \begin{equation*}
 \begin{tikzpicture}
 \node at (0,0) [name=A, anchor=east] {\ensuremath{\Pi\colon \ \widetilde{\Omega T}}};
 \node at (1,0) [name=B, anchor=west] {\ensuremath{\mP_1T},};
 \node at (0,-0.8) [name=C, anchor=east] {\ensuremath{(z,\gamma)}};
 \node at (1,-0.8) [name=D, anchor=west] {\ensuremath{\gamma,}};
 \draw[->] (A) -- (B);
 \draw[|->] (C) -- (D);
 \end{tikzpicture}
 \end{equation*}
 where $\widetilde{\Omega T}$ is the central extension with cocycle $c\inv$, equipped with the $\mP_1T$-action
 \begin{gather*}
 (z,\gamma)^\varphi=(ez,\gamma),
 \end{gather*}
 with
 \begin{gather*}
 e = \exp\(\int_0^1I(f(t),\dot g(t)) {\rm d}t\).
 \end{gather*}
 Here $(f,g)$ is the lift of $(\varphi,\gamma)$ to $\mathfrak t\times\mathfrak t$ that starts at $(0,0)$.
\end{Lemma}
\begin{proof}
 The kernel of the source of $\mG$ consists of triples $(z,1,\gamma)$ where $\gamma$ is a closed loop and $1$ denotes the constant loop. Since~$c$ is invariant under reparametrization, we have
 \begin{gather*}
 c(1*\overline\gamma,1*\overline\beta) = c(\overline\gamma,\overline\beta) = c(\gamma,\beta)\inv.
 \end{gather*}
 Similarly, the conjugation action of the loop $\varphi *\overline\varphi$ on $(z,1,\gamma)$ is calculated
 as in the proof of Lemma~\ref{lemma7.1}. 
\end{proof}

In the following, we deviate from the notation in Section~\ref{sec:Extraspecial} and write $\widetilde\Lat$ for the central extension of $\Lat$ by $U(1)$ with cocycle $(-1)^J$. Recall that a~(weak) map of crossed modules, in the sense of \cite[Def\/inition~8.4]{Noohi:Notes_on_2-groupoids}, is an equivalence if it induces a monoidal equivalence of the corresponding categorical groups.
\begin{Lemma}
 The crossed module of the previous lemma is equivalent to the crossed module
 \begin{equation*}
 \begin{tikzpicture}
 \node at (0,0) [name=A, anchor=east] {\ensuremath{\Xi\colon \ \text{\rm $\widetilde\Lat$}}};
 \node at (1,0) [name=B, anchor=west] {\ensuremath{\mathfrak t,}};
 \node at (0,-0.8) [name=C, anchor=east] {\ensuremath{(z,m)}};
 \node at (1,-0.8) [name=D, anchor=west] {\ensuremath{m}};
 \draw[->] (A) -- (B);
 \draw[|->] (C) -- (D);
 \end{tikzpicture}
 \end{equation*}
 with action
 \begin{gather*}
 (z,m)^x=\big(z\cdot\exp\big(\tfrac12I(x,m)\big),m\big).
 \end{gather*}
\end{Lemma}
\begin{proof}
 For $x\in\mathfrak t$, def\/ine the path
 \begin{equation*}
 \begin{tikzpicture}
 \node at (0,0) [name=A, anchor=east] {\ensuremath{\gamma_x\colon \ [0,1]}};
 \node at (1,0) [name=B, anchor=west] {\ensuremath{T,}};
 \node at (0,-0.8) [name=C, anchor=east] {\ensuremath{t}};
 \node at (1,-0.8) [name=D, anchor=west] {\ensuremath{\exp(tx)}};
 \draw[->] (A) -- (B);
 \draw[|->] (C) -- (D);
 \end{tikzpicture}
 \end{equation*}
from $1$ to $\exp(x)$ in $T$. The resulting group homomorphism from~$\mathfrak t$ to $\mP_1T$ identif\/ies $\Lat$ with the subgroup $\vT\sub\Omega T$, and makes $\Xi$ into an equivalent sub-crossed module of~$\Pi$.
\end{proof}
\begin{Lemma}
 We have an equivalence
 \begin{gather*}\Psi \simeq \Xi\end{gather*}
 between the crossed module of Construction~{\rm \ref{cons:strict}} and that of the previous lemma.
\end{Lemma}
\begin{proof} Let $\longmap a\Lat\{\pm1\}$ be any 1-cochain on $\Lat $ with boundary $(-1)^J$. Then we have the group homomorphism
 \begin{equation*}
 \begin{tikzpicture}
 \node at (0,0) [name=A, anchor=east] {\ensuremath{p_2\colon \ \Lat\times U(1)}};
 \node at (1,0) [name=B, anchor=west] {\ensuremath{\widetilde\Lat,}};
 \node at (0,-0.8) [name=C, anchor=east] {\ensuremath{(m,z)}};
 \node at (1,-0.8) [name=D, anchor=west] {\ensuremath{(z\cdot a(m),m).}};
 \draw[->] (A) -- (B);
 \draw[|->] (C) -- (D);
 \end{tikzpicture}
 \end{equation*}
 Together with the maps $p_1=\id_{\mathfrak t}$ and
 \begin{equation*}
 \begin{tikzpicture}
 \node at (0,0) [name=A, anchor=east] {\ensuremath{\eps\colon \ \mathfrak t\times\mathfrak t}};
 \node at (1,0) [name=B, anchor=west] {\ensuremath{\widetilde\Lat,}};
 \node at (0,-0.8) [name=C, anchor=east] {\ensuremath{(x,y)}};
 \node at (1,-0.8) [name=D, anchor=west] {\ensuremath{\(\exp\(\frac12J(x,y)\),0\),}};
 \draw[->] (A) -- (B);
 \draw[|->] (C) -- (D);
 \end{tikzpicture}
 \end{equation*}
 this def\/ines a weak map of crossed modules inducing an equivalence of the corresponding categorical groups.
\end{proof}

This concludes our argument that the four constructions of our categorical torus are equivalent.

\subsection*{Acknowledgments}
The author was supported by an Australian Research Fellowship and by ARC grant DP1095815. It is a pleasure to thank David Roberts for very helpful conversations and correspondence as well as his open referee report. The idea for Construction~\ref{cons:strict} came from a conversation with him, and I understand that he will also write about it elsewhere. I would like to thank Shan Shah for pointing out a mistake in an earlier version. Many thanks for helpful and inspiring conversations also go to Matthew Ando, Konrad Waldorf, Thomas Nikolaus, Geof\/frey Mason, Arun Ram, and Alex Ghitza. Finally I wish to thank the anonymous referee for helpful comments.

\pdfbookmark[1]{References}{ref}
\LastPageEnding


\begin{thebibliography}{99}
\footnotesize\itemsep=0pt

\bibitem{Baez:Huerta:An_invitation}
Baez J.C., Huerta J., An invitation to higher gauge theory, \href{https://doi.org/10.1007/s10714-010-1070-9}{\textit{Gen.
 Relativity Gravitation}} \textbf{43} (2011), 2335--2392, \href{https://arxiv.org/abs/1003.4485}{arXiv:1003.4485}.

\bibitem{Baez:Stevenson:Crans:Schreiber}
Baez J.C., Stevenson D., Crans A.S., Schreiber U., From loop groups to
 2-groups, \href{https://doi.org/10.4310/HHA.2007.v9.n2.a4}{\textit{Homology Homotopy Appl.}} \textbf{9} (2007), 101--135,
 \href{https://arxiv.org/abs/math.QA/0504123}{math.QA/0504123}.

\bibitem{Brylinski:Loop_spaces_characteristic_classes}
Brylinski J.-L., Loop spaces, characteristic classes and geometric quantization,
 \href{https://doi.org/10.1007/978-0-8176-4731-5}{\textit{Modern Birkh\"auser Classics}}, Birkh\"auser Boston, Inc., Boston, MA, 2008.

\bibitem{Brylinski:McLaughlin:The_geometry_of_degree-4_characteristic_classes_II}
Brylinski J.-L., McLaughlin D.A., The geometry of degree-{$4$} characteristic
 classes and of line bundles on loop spaces.~{II}, \href{https://doi.org/10.1215/S0012-7094-96-08305-2}{\textit{Duke Math.~J.}}
 \textbf{83} (1996), 105--139.

\bibitem{Carey:Johnson:Murray:Stevenson:Wang}
Carey A.L., Johnson S., Murray M.K., Stevenson D., Wang B.-L., Bundle gerbes for
 {C}hern--{S}imons and {W}ess--{Z}umino--{W}itten theories, \href{https://doi.org/10.1007/s00220-005-1376-8}{\textit{Comm.
 Math. Phys.}} \textbf{259} (2005), 577--613, \href{https://arxiv.org/abs/math.DG/0410013}{math.DG/0410013}.

\bibitem{Chern:Simons:Characteristic_forms_and_geometric_invariants}
Chern S.S., Simons J., Characteristic forms and geometric invariants,
 \href{https://doi.org/10.2307/1971013}{\textit{Ann. of Math.}} \textbf{99} (1974), 48--69.

\bibitem{Conway:Sloane:Sphere_packings}
Conway J.H., Sloane N.J.A., Sphere packings, lattices and groups,
 \href{https://doi.org/10.1007/978-1-4757-6568-7}{\textit{Grundlehren der Mathematischen Wissenschaften}}, Vol.~290, 3rd ed.,
 Springer-Verlag, New York, 1999.

\bibitem{Johnson-Freyd:Treumann}
Johnson-Freyd T., Treumann D., ${H}^4({\rm Co}_0; {\bf Z})= {\bf Z}/24$,
 \href{https://arxiv.org/abs/1707.07587}{arXiv:1707.07587}.

\bibitem{Murray:Roberts:Stevenson:Vozzo}
Murray M.K., Roberts D.M., Stevenson D., Vozzo R.F., Equivariant bundle gerbes,
 \href{https://doi.org/10.4310/ATMP.2017.v21.n4.a3}{\textit{Adv. Theor. Math. Phys.}} \textbf{21} (2017), 921--975,
 \href{https://arxiv.org/abs/1506.07931}{arXiv:1506.07931}.

\bibitem{Noohi:Notes_on_2-groupoids}
Noohi B., Notes on 2-groupoids, 2-groups and crossed modules, \href{https://doi.org/10.4310/HHA.2007.v9.n1.a3}{\textit{Homology
 Homotopy Appl.}} \textbf{9} (2007), 75--106, \href{https://arxiv.org/abs/math.CT/0512106}{math.CT/0512106}.

\bibitem{Pressley:Segal}
Pressley A., Segal G., Loop groups, \textit{Oxford Mathematical Monographs, Oxford
 Science Publications}, The Clarendon Press, Oxford University Press, New York,
 1986.

\bibitem{Schommer-Pries:Central_extensions}
Schommer-Pries C.J., Central extensions of smooth 2-groups and a
 f\/inite-dimensional string 2-group, \href{https://doi.org/10.2140/gt.2011.15.609}{\textit{Geom. Topol.}} \textbf{15} (2011),
 609--676, \href{https://arxiv.org/abs/0911.2483}{arXiv:0911.2483}.

\bibitem{Tits:Le_Monstre}
Tits J., Le {M}onstre (d'apr\`es {R}.~{G}riess, {B}.~{F}ischer et al.),
 \textit{Ast\'erisque} \textbf{121--122} (1985), 105--122.

\bibitem{Wagemann:Wockel:A_cocycle_model}
Wagemann F., Wockel C., A cocycle model for topological and {L}ie group
 cohomology, \href{https://doi.org/10.1090/S0002-9947-2014-06107-2}{\textit{Trans. Amer. Math. Soc.}} \textbf{367} (2015), 1871--1909,
 \href{https://arxiv.org/abs/1110.3304}{arXiv:1110.3304}.

\bibitem{Waldorf:Multiplicative_bundle_gerbes_with_connection}
Waldorf K., Multiplicative bundle gerbes with connection, \href{https://doi.org/10.1016/j.difgeo.2009.10.006}{\textit{Differential
 Geom. Appl.}} \textbf{28} (2010), 313--340, \href{https://arxiv.org/abs/0804.4835}{arXiv:0804.4835}.

\bibitem{Waldorf:A_construction_of_string_2-group_models}
Waldorf K., A construction of string 2-group models using a
 transgression-regression technique, in Analysis, Geometry and Quantum Field
 Theory, \href{https://doi.org/10.1090/conm/584/11588}{\textit{Contemp. Math.}}, Vol.~584, Amer. Math. Soc., Providence, RI,
 2012, 99--115, \href{https://arxiv.org/abs/1201.5052}{arXiv:1201.5052}.

\bibitem{Waldorf:Transgression_to_loop_spaces_and_its_inverse_II}
Waldorf K., Transgression to loop spaces and its inverse, {II}:~{G}erbes and
 fusion bundles with connection, \href{https://doi.org/10.4310/AJM.2016.v20.n1.a4}{\textit{Asian~J. Math.}} \textbf{20} (2016),
 59--115, \href{https://arxiv.org/abs/1004.0031}{arXiv:1004.0031}.

\end{thebibliography}
\end{document}